\documentclass{amsart}

\usepackage{amsfonts,amssymb,verbatim,amsmath,amsthm,latexsym,textcomp,amscd, hyperref}
\usepackage{latexsym,amsfonts,amssymb,epsfig,verbatim}
\usepackage{amsmath,amsthm,amssymb,latexsym,graphics,textcomp}
\usepackage{graphicx}
\usepackage{color}
\usepackage{url}

\input xy
\xyoption{all}

\newcommand{\Appendix}[1]{%
  \refstepcounter{section}%
  \addcontentsline{toc}{section}%
    {\bfseries\appendixname~\thesection\ #1}%
    {\medskip\noindent \Large\bfseries\appendixname\ \thesection\ #1}%
\sectionmark{#1}\smallskip\noindent
\renewcommand{\theequation}{{\bf 
{{\thesection}.\arabic{subsection}}.{\arabic{equation}}}}
}
\begin{document}

\newtheorem{theorem}{Theorem}[section]
\newtheorem{prop}[theorem]{Proposition}
\newtheorem{lemma}[theorem]{Lemma}
\newtheorem{cor}[theorem]{Corollary}
\newtheorem{definition}[theorem]{Definition}
\newtheorem{defn}[theorem]{Definition}
\newtheorem{conj}[theorem]{Conjecture}
\newtheorem{claim}[theorem]{Claim}
\newtheorem{remark}[theorem]{Remark}
\newtheorem{rmk}[theorem]{Remark}
\newtheorem{defth}[theorem]{Definition-Theorem}

\newcommand{\boundary}{\partial}
\newcommand{\C}{{\mathbb C}}
\newcommand{\integers}{{\mathbb Z}}
\newcommand{\natls}{{\mathbb N}}
\newcommand{\ratls}{{\mathbb Q}}
\newcommand{\reals}{{\mathbb R}}
\newcommand{\proj}{{\mathbb P}}
\newcommand{\lhp}{{\mathbb L}}
\newcommand{\tube}{{\mathbb T}}
\newcommand{\cusp}{{\mathbb P}}
\newcommand\AAA{{\mathcal A}}
\newcommand\BB{{\mathcal B}}
\newcommand\CC{{\mathcal C}}
\newcommand\DD{{\mathcal D}}
\newcommand\EE{{\mathcal E}}
\newcommand\FF{{\mathcal F}}
\newcommand\GG{{\mathcal G}}
\newcommand\HH{{\mathcal H}}
\newcommand\II{{\mathcal I}}
\newcommand\JJ{{\mathcal J}}
\newcommand\KK{{\mathcal K}}
\newcommand\LL{{\mathcal L}}
\newcommand\MM{{\mathcal M}}
\newcommand\NN{{\mathcal N}}
\newcommand\OO{{\mathcal O}}
\newcommand\PP{{\mathcal P}}
\newcommand\QQ{{\mathcal Q}}
\newcommand\RR{{\mathcal R}}
\newcommand\SSS{{\mathcal S}}
\newcommand\TT{{\mathcal T}}
\newcommand\UU{{\mathcal U}}
\newcommand\VV{{\mathcal V}}
\newcommand\WW{{\mathcal W}}
\newcommand\XX{{\mathcal X}}
\newcommand\YY{{\mathcal Y}}
\newcommand\ZZ{{\mathcal Z}}
\newcommand\CH{{\CC\HH}}
\newcommand\MF{{\MM\FF}}
\newcommand\PEY{{\PP\EE\YY}}
\newcommand\RCT{{{\mathcal R}_{CT}}}
\newcommand\PMF{{\PP\kern-2pt\MM\FF}}
\newcommand\FL{{\FF\LL}}
\newcommand\PML{{\PP\kern-2pt\MM\LL}}
\newcommand\GL{{\GG\LL}}
\newcommand\Pol{{\mathcal P}}
\newcommand\half{{\textstyle{\frac12}}}
\newcommand\Half{{\frac12}}
\newcommand\Mod{\operatorname{Mod}}
\newcommand\Area{\operatorname{Area}}
\newcommand\ep{\epsilon}
\newcommand\hhat{\widehat}
\newcommand\Proj{{\mathbf P}}
\newcommand\U{{\mathbf U}}
 \newcommand\Hyp{{\mathbf H}}
\newcommand\D{{\mathbf D}}
\newcommand\Z{{\mathbb Z}}
\newcommand\R{{\mathbb R}}
\newcommand\Q{{\mathbb Q}}
\newcommand\E{{\mathbb E}}
\newcommand\til{\widetilde}
\newcommand\length{\operatorname{length}}
\newcommand\tr{\operatorname{tr}}
\newcommand\gesim{\succ}
\newcommand\lesim{\prec}
\newcommand\simle{\lesim}
\newcommand\simge{\gesim}
\newcommand{\simmult}{\asymp}
\newcommand{\simadd}{\mathrel{\overset{\text{\tiny $+$}}{\sim}}}
\newcommand{\ssm}{\setminus}
\newcommand{\diam}{\operatorname{diam}}
\newcommand{\pair}[1]{\langle #1\rangle}
\newcommand{\T}{{\mathbf T}}
\newcommand{\inj}{\operatorname{inj}}
\newcommand{\pleat}{\operatorname{\mathbf{pleat}}}
\newcommand{\short}{\operatorname{\mathbf{short}}}
\newcommand{\vertices}{\operatorname{vert}}
\newcommand{\collar}{\operatorname{\mathbf{collar}}}
\newcommand{\bcollar}{\operatorname{\overline{\mathbf{collar}}}}
\newcommand{\I}{{\mathbf I}}
\newcommand{\tprec}{\prec_t}
\newcommand{\fprec}{\prec_f}
\newcommand{\bprec}{\prec_b}
\newcommand{\pprec}{\prec_p}
\newcommand{\ppreceq}{\preceq_p}
\newcommand{\sprec}{\prec_s}
\newcommand{\cpreceq}{\preceq_c}
\newcommand{\cprec}{\prec_c}
\newcommand{\topprec}{\prec_{\rm top}}
\newcommand{\Topprec}{\prec_{\rm TOP}}
\newcommand{\fsub}{\mathrel{\scriptstyle\searrow}}
\newcommand{\bsub}{\mathrel{\scriptstyle\swarrow}}
\newcommand{\fsubd}{\mathrel{{\scriptstyle\searrow}\kern-1ex^d\kern0.5ex}}
\newcommand{\bsubd}{\mathrel{{\scriptstyle\swarrow}\kern-1.6ex^d\kern0.8ex}}
\newcommand{\fsubeq}{\mathrel{\raise-.7ex\hbox{$\overset{\searrow}{=}$}}}
\newcommand{\bsubeq}{\mathrel{\raise-.7ex\hbox{$\overset{\swarrow}{=}$}}}
\newcommand{\tw}{\operatorname{tw}}
\newcommand{\base}{\operatorname{base}}
\newcommand{\trans}{\operatorname{trans}}
\newcommand{\rest}{|_}
\newcommand{\bbar}{\overline}
\newcommand{\UML}{\operatorname{\UU\MM\LL}}
\newcommand{\EL}{\mathcal{EL}}
\newcommand{\tsum}{\sideset{}{'}\sum}
\newcommand{\tsh}[1]{\left\{\kern-.9ex\left\{#1\right\}\kern-.9ex\right\}}
\newcommand{\Tsh}[2]{\tsh{#2}_{#1}}
\newcommand{\qeq}{\mathrel{\approx}}
\newcommand{\Qeq}[1]{\mathrel{\approx_{#1}}}
\newcommand{\qle}{\lesssim}
\newcommand{\Qle}[1]{\mathrel{\lesssim_{#1}}}
\newcommand{\simp}{\operatorname{simp}}
\newcommand{\vsucc}{\operatorname{succ}}
\newcommand{\vpred}{\operatorname{pred}}
\newcommand\fhalf[1]{\overrightarrow {#1}}
\newcommand\bhalf[1]{\overleftarrow {#1}}
\newcommand\sleft{_{\text{left}}}
\newcommand\sright{_{\text{right}}}
\newcommand\sbtop{_{\text{top}}}
\newcommand\sbot{_{\text{bot}}}
\newcommand\sll{_{\mathbf l}}
\newcommand\srr{_{\mathbf r}}
\newcommand\geod{\operatorname{\mathbf g}}
\newcommand\mtorus[1]{\boundary U(#1)}
\newcommand\A{\mathbf A}
\newcommand\Aleft[1]{\A\sleft(#1)}
\newcommand\Aright[1]{\A\sright(#1)}
\newcommand\Atop[1]{\A\sbtop(#1)}
\newcommand\Abot[1]{\A\sbot(#1)}
\newcommand\boundvert{{\boundary_{||}}}
\newcommand\storus[1]{U(#1)}
\newcommand\Momega{\omega_M}
\newcommand\nomega{\omega_\nu}
\newcommand\twist{\operatorname{tw}}
\newcommand\modl{M_\nu}
\newcommand\MT{{\mathbb T}}
\newcommand\Teich{{\mathcal T}}
\renewcommand{\Re}{\operatorname{Re}}
\renewcommand{\Im}{\operatorname{Im}}

\title{Ending Laminations and Cannon-Thurston Maps }

\author[Mahan Mj with Shubhabrata Das]{Mahan Mj \\ with an Appendix by Shubhabrata Das and Mahan Mj}

\address{RKM Vivekananda University, Belur Math, WB-711 202, India}

\email{mahan.mj@gmail.com; mahan@rkmvu.ac.in}
\email{shubhabrata.gt@gmail.com}

\subjclass[2010]{57M50, 20F67 (Primary); 20F65,  22E40  (Secondary)}

\begin{abstract}

In earlier work, we had shown that Cannon-Thurston
 maps exist for Kleinian  surface groups without accidental parabolics. 
In this paper we prove that pre-images of points are precisely end-points of leaves of the ending
lamination whenever the Cannon-Thurston map is not one-to-one.  
\end{abstract}

\maketitle

\tableofcontents

\section{Introduction}

\subsection{Statement of Results}
In earlier work  we showed:

\begin{theorem} \cite{mahan-split}
 Let
$\rho : \pi_1(S) \rightarrow
PSL_2(C)$ be a discrete faithful representation of a
surface group with or without punctures,
and without accidental parabolics. Let $M =
{{\mathbb{H}}^3}/{\rho (\pi_1 (S))}$.
Let $i$ be an embedding of $S$
in $M$ that induces a homotopy equivalence. Then the embedding
$\tilde{i} : \widetilde{S} \rightarrow \widetilde{M} = {\mathbb{H}}^3$
extends continuously to a map $\hat{i}: {\mathbb{D}}^2 \rightarrow
{\mathbb{D}}^3$. Further, the limit set of ${\rho (\pi_1 (S))}$ is
locally connected.
\label{split}
\end{theorem}

This generalizes the first part of the next theorem  due to Cannon and Thurston \cite{CT} (for 3 manifolds fibering over the circle) and Minsky
 \cite{minsky-jams} (for bounded geometry closed surface Kleinian groups): \\
\begin{theorem} \cite{CT, CTpub, minsky-jams}
Suppose a closed surface group $\pi_1 (S)$ of bounded geometry acts freely and properly
discontinuously on ${\mathbb{H}}^3$ by hyperbolic isometries. Then the inclusion
$\tilde{i} : \widetilde{S} \rightarrow {\mathbb{H}}^3$ extends
continuously to the boundary. Further, pre-images of points on the boundary  are precisely ideal boundary points of a leaf of the ending
lamination, or ideal boundary points of
a complementary ideal polygon whenever the Cannon-Thurston map is not one-to-one.
\label{ct}
\end{theorem}

In the main body of this paper, we generalize the second
part of the
above theorem to arbitrary Kleinian closed surface groups without accidental parabolics.

\begin{theorem}
Suppose a closed surface group $\pi_1 (S)$  acts freely and properly
discontinuously on ${\mathbb{H}}^3$ by hyperbolic isometries. Then the inclusion
$\tilde{i} : \widetilde{S} \rightarrow {\mathbb{H}}^3$ extends
continuously to the boundary. Further, pre-images of points
on the boundary  are precisely ideal boundary points of a leaf of the ending
lamination, or ideal boundary points of
a complementary ideal polygon whenever the Cannon-Thurston map is not one-to-one.
\label{main}
\end{theorem}

In passing from Theorem \ref{ct} to Theorem \ref{main} we have  removed the hypothesis of {\em bounded geometry}. In an appendix to the paper
we extend Theorem \ref{main} to the case of  surfaces with cusps (Theorem \ref{main-cusp}).

\subsection{Outline and Applications}\label{out} We first outline the main steps involved in the proof of the main Theorem \ref{main}. To fix notions, we let $M$ 
be the convex core of a simply or doubly degenerate hyperbolic 3-manifold homotopy equivalent to a surface $S$. We also assume that
an inclusion $i: S \rightarrow M$ inducing the homotopy equivalence is fixed. Let $\til{i} : \til{S} \rightarrow \til{M}$ denote the lift of $i$ to universal covers.

\medskip

\noindent {\bf Recapitulation of Theorem \ref{split} from
\cite{mahan-split}:} \\ 
To show that a Cannon-Thurston map exists we have to show that $\til i$ extends continuously to the boundary giving $\hat{i}: {\mathbb{D}}^2 \rightarrow
{\mathbb{D}}^3$.
 The proof of the main Theorem \ref{split} of
\cite{mahan-split} proceeds (cf. Lemma \ref{crit} below) by showing that given a geodesic segment $\lambda$ in (the intrinsic metric on)
 $\widetilde{S}$ lying outside a large ball about a fixed reference  point $o$ in $\til S$, the hyperbolic geodesic in 
$\til M$ joining its end-points lies outside a large ball
about $\til{i} (o)$ in $\til M$. Towards this a {\em hyperbolic ladder} ${\mathcal{L}}_\lambda$ is constructed in $\til M$ containing $\lambda$ satisfying 
the following:\\
a) a 
(weak) quasiconvexity property, \\
b) If $\lambda$ lies outside a large ball about $o$ in the intrinsic metric on
 $\widetilde{S}$, then ${\mathcal{L}}_\lambda$ lies outside a large ball about $\til{i}(o)$ in $\til M$. \\

The quasiconvexity property of ${\mathcal{L}}_\lambda$ ensures control over the hyperbolic geodesic in 
$\til M$ joining the end-points of $\til{i} (\lambda)$. In particular, if ${\mathcal{L}}_\lambda$ lies outside a large ball about $\til{i}(o)$ in $\til M$
then so does the geodesic in $\til M$ joining the end-points of $\til{i} (\lambda)$. This guarantees the existence of the Cannon-Thurston map
$\hat{i}$ in Theorem \ref{split}. 

\medskip

\noindent {\bf Scheme of proof of Theorem \ref{main}:} \\ Theorem \ref{main} builds on Theorem \ref{split} by describing the structure of the
Cannon-Thurston map obtained in \cite{mahan-split}. The crux of the proof of Theorem \ref{main} involves an analysis of the structure of certain
specific ladders ${\mathcal{L}}_\lambda$. The existence and weak quasiconvexity of these ladders was shown in \cite{mahan-split}, but
the analysis (see Steps 2, 3 below) was missing. In fact, even for punctured torus groups, where the existence of Cannon-Thurston maps was shown by
McMullen \cite{ctm-locconn}, Theorem \ref{main} is new.
We now proceed with a step-by-step outline of the  proof of Theorem \ref{main}. 

Step 1) The easy part (Section \ref{easy}) of  Theorem \ref{main} consists in showing that the end-points of leaves of ending laminations are identified by the Cannon-Thurston map
$\hat{i}$. The essential point is that a leaf of the ending lamination in $\til S$ can be approximated by the lifts to $\til S$
 of a sequence of  closed  curves in $S$ whose geodesic realizations exit the relevant end of $M$. We shall refer to this step as the {\bf forward direction}
of Theorem  \ref{main}.

Step 2) The hard part of the proof  (Section \ref{closed} and Appendix \ref{appendix}) consists in showing that if the Cannon-Thurston map
$\hat i$ identifies a pair of points, then they are the ideal end-points of a leaf of the ending lamination or an ideal complementary polygon.
We shall refer to this step as the {\bf reverse direction}
of Theorem  \ref{main}.
Bi-infinite geodesics whose end-points are identified by $\hat i$ are referred to as CT leaves (cf. Section \ref{lamns}). 

The proof proceeds by analyzing the structure of the ladder ${\mathcal{L}}_\lambda$  for $\lambda$ a CT leaf. The heart of the proof lies in
Proposition \ref{qgeodasymp} (Asymptotic Quasigeodesic Rays) which essentially says that "vertical" quasigeodesic rays lying on such a ladder 
 ${\mathcal{L}}_\lambda$  are all asymptotic to some point $z_\lambda \in \partial {\mathbb{H}}^3$. Further the end-points of $\lambda$ are identified with $z_\lambda$
under $\hat i$.

Step 3) Given Proposition \ref{qgeodasymp} there are two ways to complete the proof of Theorem \ref{main}: \\
a) Look at the action of $\pi_1(S)$ on the $\mathbb{R}$-tree dual to the ending lamination. If there is a CT leaf that is not a leaf of the ending lamination,
then we construct a CT leaf  (Section \ref{closed}) whose ideal end-points consist of the attracting and repelling fixed points
$g^{-\infty} , g^\infty$ for some $g \in \pi_1(S)$. This is a contradiction as $g$ is a hyperbolic (loxodromic) element. This is the approach taken in
Section \ref{closed}. \\
b) Alternately use Proposition \ref{qgeodasymp} and a Lemma of Bowditch (Lemma 9.2 of \cite{bowditch-ct})
to show that the collection of CT-leaves forms a lamination. Since the easy direction shows that the ending lamination is contained in the collection of CT leaves, 
this forces the collection of CT-leaves to exactly equal the ending lamination.  This is the approach taken in Appendix \ref{appendix}.

\medskip

\noindent {\bf Applications:} \\
1) We prove the following strengthening of a rigidity Theorem due to Brock-Canary-Minsky \cite{minsky-elc2}:

{\bf Theorem \ref{rig}:} {\it
Let $G$ be a  closed surface group.
Let $\rho (G) = \Gamma$ and $\rho_1 (G) = \Gamma_1$ be two  simply or doubly degenerate  representations of
  $G$ into $PSl_2(\mathbb{C})$ with limits sets $\Lambda, \Lambda_1$. Suppose that the $G-$ actions on $\Lambda, \Lambda_1$ are topologically conjugate.
Then   $\rho$ and $\rho_1$ are
  quasiconformally conjugate.}

\smallskip

2) Theorem \ref{main} and its generalization Theorem \ref{main-cusp} are used to prove discreteness
of commensurators of finitely generated infinite covolume Kleinian groups in \cite{llr} and \cite{mahan-discrete}.\\
3) Theorems \ref{main} and  \ref{main-cusp} are extended to arbitrary finitely generated Kleinian groups in \cite{mahan-kl}.

\medskip

\noindent {\bf Organization of the paper:} \\ Section \ref{prelim} of the paper deals with preliminary concepts and material from \cite{mahan-split}.
Section \ref{easy} proves the easy direction of Theorems \ref{main} and \ref{main-cusp}: End-points of leaves of the ending
lamination are identified by the Cannon-Thurston map. The arguments in Sections \ref{prelim} and \ref{lamns} give a unified treatment 
for surfaces with or without cusps. Section \ref{closed} proves the harder direction of Theorem \ref{main} for surfaces without cusps. A slight modification
of a fact proven for closed surfaces (Remark \ref{twin-cusp}) will be used for cusped surfaces. We indicate this in Section \ref{closed} itself.
Appendix \ref{appendix}
 deals
with surfaces with cusps.

\medskip

\noindent {\bf Acknowledgments:} The author would like to thank
the referee(s) 
for pointing out  errors and omissions and for suggesting corrections. The research for the case without parabolics
is supported in part by a DST research grant
DyNo. 100/IFD/8347/2008-2009.  The research for the case with parabolics
is supported in part by a CEFIPRA project grant
4301-1.

\section{Preliminaries}\label{prelim}

\subsection{Hyperbolic Metric Spaces}\label{hypsec}

Let $(X,{d_X})$ be a hyperbolic metric space and $Y$ be a subspace that is
hyperbolic with the inherited path metric $d_Y$.
By
adjoining the Gromov boundaries $\partial{X}$ and $\partial{Y}$
 to $X$ and $Y$, one obtains their compactifications
$\widehat{X}$ and $\widehat{Y}$ respectively.

Let $ i :Y \rightarrow X$ denote inclusion.

\begin{definition}   Let $X$ and $Y$ be hyperbolic metric spaces and
$i : Y \rightarrow X$ be an embedding.
 A {\bf Cannon-Thurston map} $\hat{i}$  from $\widehat{Y}$ to
 $\widehat{X}$ is a continuous extension of $i$.
\end{definition}

The following  lemma (Lemma 2.1 of \cite{mitra-ct})
 says that a Cannon-Thurston map exists
if and only if for all $M > 0$ and $y \in Y$, there exists $N > 0$ such that if a geodesic $\lambda$ in $Y$
lies outside an $N$ ball around $y$ in $Y$,  then
any geodesic in $X$ joining the end-points of $\lambda$ lies
outside the $M$ ball around $i(y)$ in $X$. An equivalent statement is that the Cannon-Thurston map exists if and only if sets of small visual diameter go
to sets of small visual diameter.

\begin{lemma}\label{crit} (Lemma 2.1 of \cite{mitra-ct}) Let $i: Y \rightarrow X$ be an inclusion of hyperbolic metric spaces. A Cannon-Thurston map from $\widehat{Y}$ to $\widehat{X}$
 exists if and only if the following condition is satisfied:\\
Given ${y_0}\in{Y}$, there exists a non-negative function  $M(N)$, such that 
 $M(N)\rightarrow\infty$ as $N\rightarrow\infty$, and such that for all geodesic segments
 $\lambda$  lying outside an $N$-ball
around ${y_0}\in{Y}$,   any geodesic segment in $X$ joining
the end-points of $i(\lambda)$ lies outside the $M(N)$-ball around 
$i({y_0})\in{X}$.
\end{lemma}

\noindent {\bf Relative Hyperbolicity and Electric Geometry}
We refer the reader to  \cite{farb-relhyp}
for terminology and details on relative
hyperbolicity and electric geometry. \\
Let $X$ be a $\delta$-hyperbolic metric
space, and $\mathcal{H}$ a family of $C$-quasiconvex, $D$-separated,
 collection of subsets. 
Recall
\cite{mahan-split} that {\bf electrocution} of the collection $\mathcal{H}$ in $X$ means constructing
an auxiliary space $X_{el} = X \bigcup_{H \in \mathcal{H}} (H \times I)$
with $H \times \{ 0 \}$ identified to $H \subset X$ and  
$H \times \{ 1 \}$ equipped with the zero metric. This is a geometric
`coning' construction.

Then by work of Farb \cite{farb-relhyp},
$X_{el}$ obtained by electrocuting the subsets in $\mathcal{H}$ is
a $\Delta = \Delta ( \delta , C, D)$ -hyperbolic metric space. 

 Now,
let $\alpha = [a,b]$ be a hyperbolic geodesic in $X$ and $\beta $ be
an electric
$P$-quasigeodesic without backtracking
 joining $a, b$. Starting from the left
 of
$\beta$, replace each maximal subsegment, (with end-points $p,
 q$, say)
 lying within some $H\times \{ 1\} (H \in \mathcal{H})$
by a hyperbolic  geodesic $[p,q]$. The resulting
{\bf connected}
path $\beta_q$ is called an {\em electro-ambient representative} in
$X$. Electro-ambient representatives are useful in light of the following.

 \begin{lemma} \cite{mahan-split}\label{ea-genl} Given $\delta, C \geq 0, D >0$, there exists $D_0$ such that the following holds: \\
Let $(X,d)$ be a $\delta$-hyperbolic metric
space, and $\mathcal{H}$ a family of $C$-quasiconvex, $D$-separated,
 collection of subsets. 
Let $(X_{el},d_{el})$ denote the electric space obtained from $X$ by electrocuting the family 
$\mathcal{H}$. Let $\alpha = [a,b]$ be a  geodesic in $X$ and $\beta_q $ be
an electro-ambient representative of an electric geodesic joining $a, b$. Then $\alpha , \beta_q $
lies within a bounded distance $D_0$ of each other in $(X_{el},d_{el})$. Further, $\alpha$ lies within a bounded distance $D_0$ of $ \beta_q $
in $(X,d)$.
\end{lemma}

\noindent {\bf  Partial Electrocution}\\
Let $Y$ be the convex core of a simply (resp. doubly) degenerate 3-manifold with cusps. 
After removing an open neighborhood  of cusps we get a manifold with boundary of the form $S \times J$, where $J = [0, \infty )$ (resp.
$\mathbb R$), where $S$ is a surface with boundary. Surfaces minus open neighborhoods of cusps shall sometimes be referred to as
{\it truncated surfaces}.
Let $\mathcal{B}$ denote the equivariant collection of horoballs in $\til Y$ covering the cusps of $Y$. Let $X$ denote
$\til Y$ minus the interior of the horoballs in $\mathcal{B}$. Let 
$\mathcal{H}$ denote the collection of boundary horospheres. Then each
$H \in \mathcal{H}$ with the induced metric is isometric to a Euclidean
product $\mathbb{R} \times J$. 
 We shall need to  equip  each 
$H \in \HH$ is with a new pseudo-metric called the {\it partially electrocuted metric}
 by giving it the product of the zero metric (in the $\mathbb R$-direction) with the Euclidean metric
 (in the $J$-direction). The resulting space is quasi-isometric to what one would get by gluing
to each $H$ the mapping cylinder of the projection of $H$ onto the $J$-factor. Let $\JJ$ denote the collection of copies of $J$ obtained in this 
construction and let $(\PEY , d_{pel})$ denote the resulting partially electrocuted space. (See \cite{mahan-pal} for a more general discussion.)
We have the following basic Lemma.

\begin{lemma} \cite{mahan-pal}
$(\PEY ,d_{pel})$ is a hyperbolic metric space and the sets $J_\alpha \in \JJ$
are uniformly quasiconvex.
\label{pel}
\end{lemma}

\subsection{Split Geometry}

We shall briefly recall the essential aspects of split geometry from \cite{minsky-elc1, mahan-split}.
 We shall also need the construction of certain quasiconvex ladder-like sets ${\LL}_\lambda$. Since we shall deal with surfaces with cusps (or punctures) in the Appendix, we give a unified
exposition for surfaces with or without punctures. If a finite area hyperbolic
surface has cusps, we shall remove an open neighborhood of the cusp and denote the resulting {\it truncated surface} by $S$.
  In this subsection therefore $S$ will denote a compact surface, possibly with boundary.

\smallskip

\noindent {\bf Split level Surfaces}\\
A   {\bf pants decomposition} of a compact surface $S$, possibly with boundary,
is a disjoint collection of
3-holed spheres $P_1, \cdots , P_n$ embedded in $S$ such that $S \setminus \bigcup_i P_i$ is a disjoint collection of non-peripheral
annuli in $S$, no two of which are  homotopic. 

\begin{comment}
We shall conflate a pants decomposition of $S$ with the  collection of
(isotopy classes of) {\it non-peripheral}  boundary curves of  $P_1, \cdots , P_n$. Thus when we refer to a pair of pants in a pants decomposition 
 $P_1, \cdots , P_n$ of $S$ we are referring to one of the $P_i$'s, and when we refer to a curve in a pants decomposition 
 of $S$ we are referring to one of the  non-peripheral  boundary curves of one of the $P_i$'s.
 \end{comment}

Let $N$ be the convex core of a  hyperbolic 3-manifold minus an open neighborhood of the cusp(s). Then any end $E$ of $N$ is simply degenerate
 \cite{agol-tameness, gab-cal, canary} and
 homeomorphic to $S \times [0, \infty )$, where $S$ is a compact surface, possibly with boundary. A closed geodesic in an end $E$ 
 homeomorphic to $S \times [0, \infty )$  is {\bf unknotted} if it is isotopic in $E$ to a simple closed curve in $S \times \{0 \}$ via the homeomorphism.
A {\bf tube} in an end $E \subset N$ is a regular $R-$neighborhood
$N(\gamma, R)$ of an unknotted geodesic $\gamma$ in $E$.

\begin{comment}
Let $ \theta , \omega $ be positive real numbers. 
A neighborhood $N_\epsilon (\gamma )$ of a closed geodesic $\gamma (\subset N)$
is called a $(\theta, \omega)$-thin tube if the length of $\gamma$ is less than $\theta$ and the length of the shortest
geodesic on $\partial N_\epsilon (\gamma )$ is greater than $\omega$.
\end{comment}

Let $\TT$ denote a collection of disjoint, uniformly separated tubes in ends of  $N$ 
such that\\
a) all Margulis tubes in $E$ belong to $\TT$ for all ends $E$ of $N$.\\
b) there exists $\epsilon_0 >0$ such that the injectivity radius $injrad_x(E) > \epsilon_0$ for all $x \in E \setminus  \bigcup_{T \in \TT} Int(T)$ and all ends $E$ of $N$.\\

Let $F: N \rightarrow M$ be a  bi-Lipschitz homeomorphism and let $M(0)$
be the image of $N \setminus \bigcup_{T \in \TT} Int(T)$ in $M$ under  the 
 bi-Lipschitz homeomorphism $F$. Let $\partial M(0)$ (resp. $\partial M$) denote the  boundary  of $M(0)$ (resp. $M$). Following \cite{mahan-split},
$M$ will be called the
model manifold. The metrics on $M$ and $\til M$ will be denoted by $d_M$. In \cite{minsky-elc1, minsky-elc2}, the model manifold refers to $M$ with
considerable additional structure. In particular it involves the decomposition of $M(0)$ into pieces of the form $S_{0,4} \times I$ and 
 $S_{1,1} \times I$ where $S_{0,4}$ and  $S_{1,1}$ refer to a sphere with 4 holes and a torus with one hole respectively. To distinguish between the 
 model manifold in  \cite{minsky-elc1, minsky-elc2} and that in this paper we shall refer to the former as the Minsky model. It should be pointed out that
Minsky model was constructed by Minsky in   \cite{minsky-elc1} and proven to be bi-Lipschitz homeomorphic to the hyperbolic manifold $N$ by Brock-Canary-Minsky
in    \cite{minsky-elc2}.

Let $(Q, \partial Q)$ be the unique hyperbolic  pair of pants such that each  component
of $\partial Q$ has length one. $Q$ will be called
the {\it standard} pair of pants.
An isometrically embedded copy of $(Q, \partial Q)$ in $(M(0), \partial M(0))$ will be said to be {\it flat}.

\begin{comment}
The product of the unit circle with the unit interval, 
 $S^1 \times [0,1]$ will be called the {\it standard annulus}.
\end{comment}

\begin{defn} {\rm A {\bf split level surface} associated to a pants decomposition $\{ Q_1, \cdots , Q_n \}$ of $S$ in $M(0) \subset M$
is an embedding $f : \cup_i (Q_i, \partial Q_i) \rightarrow (M(0), \partial M(0))$ such that \\
1) Each $f (Q_i, \partial Q_i)$ is flat \\
2) $f$ extends to an embedding (also denoted $f$) of $S$ into $M$ such that the interior of each annulus component of
$f(S \setminus \bigcup_i Q_i)$  lies entirely in $F(\bigcup_{T \in \TT} Int(T))$. \\
} \end{defn}

Let
$S_{i}^{s}$ denote the union of the collection of flat pairs of pants
 in the image of the embedding  $f_{i}$.

The class   of {\it all} topological embeddings from $S$ to $M$ that agree with a split level surface $f$ 
associated to a pants decomposition $\{ Q_1, \cdots , Q_n \}$ on 
$Q_1 \cup \cdots \cup Q_n$ will be denoted by $[f]$. 

We define a partial order $\leq_E$ on the collection of split level surfaces in an end $E$ of $M$ as follows: \\
$f_1 \leq_E f_2$ if there exist $g_i \in [f_i]$, $i=1,2$, such that $g_2(S)$ lies in the unbounded component of $E \setminus g_1(S)$.

A sequence $f_i$ of split level surfaces is said to exit an end $E$ if $i<j$ implies $f_i \leq_E f_j$ and further for all compact subsets $B \subset E$, there exists
$L>0$ such that $f_i(S) \cap B = \emptyset$ for all $i \geq L$.

\begin{definition}  A curve $v$ in $S \subset E$ is {\bf $l$-thin} if the core curve of the Margulis tube $T_v (\subset E \subset N)$ has length less than or equal to $l$. 
A   tube $T\in \TT$  is   $l$-thin  if its core curve    is   $l$-thin.  A   tube $T\in \TT$  is   $l$-thick if it is not    $l$-thin.  \\
A curve $v$ is said to split a pair of split level surfaces $S_i$ and $S_j$ ($i<j$) if $v$ occurs as a boundary curve of
 both $S_i$ and $S_{j}$.\\
\end{definition}

The collection of all  $l$-thin tubes is denoted as $\TT_l$. The union of all  $l$-thick tubes along with $M(0)$  is denoted as $M(l)$.

\begin{defn}
A pair of split level surfaces $S_i$ and $S_j$ ($i<j$) is said to be {\bf $k$-separated} if \\
a) for all $x \in S_i^s$, 
$d_M(x,S_j^s) \geq k$\\
b)Similarly, for all $x \in S_j^s$, $d_M(x,S_i^s) \geq k$. \end{defn}

\begin{defn} {\rm An $L$-bi-Lipschitz {\bf split  surface} in $M(l)$ associated to a pants decomposition $\{ Q_1, \cdots , Q_n \}$ of $S$
and a collection $\{ A_1, \cdots , A_m \}$ of complementary annuli (not necessarily all of them) in $S$ 
is an embedding $f : \cup_i Q_i \bigcup  \cup_i A_i \rightarrow M(l)$ such that\\
1) the restriction  $f: \cup_i (Q_i, \partial Q_i) \rightarrow (M(0), \partial M(0))$ is a split level surface \\
2) the restriction $f: A_i \rightarrow M(l)$ is an $L$-bi-Lipschitz embedding.\\
3)  $f$ extends to an embedding (also denoted $f$) of $S$ into $M$ such that the interior of each annulus component of
$f(S \setminus (\cup_i Q_i \bigcup  \cup_i A_i))$  lies entirely in $F(\bigcup_{T \in \TT_l} Int(T))$.}\end{defn}

\noindent {\bf Note:} The difference between a split level surface and a split surface is that the latter may contain
bi-Lipschitz annuli in addition to flat pairs of pants.

\smallskip

We denote split surfaces by  $\Sigma_{i}$ to distinguish them from split level surfaces $S_i$.
Let
$\Sigma_{i}^{s}$ denote the union of the collection of flat pairs of pants
and bi-Lipschitz annuli in the image of the split surface (embedding)  $\Sigma_{i}$.

\begin{theorem} \cite[Theorem 4.8]{mahan-split}
Let $N, M, M(0), S, F$ be as above and $E$ an end of $M$. For any $l$ less than the Margulis constant,
let $M(l) = \{ F(x) : {\rm injrad_x} (N) \geq l \}$. Fix a hyperbolic metric on $S$ such that each component of $\partial S$ is 
totally geodesic of length one.
 There exist $ L_1 \geq 1$, $  \epsilon_1 > 0$, $n \in \natls$, 
 and a sequence $\Sigma_i$ of $L_1$-bi-Lipschitz, $  \epsilon_1$-separated split  surfaces exiting the end $E$ of $M$
such that for all $i$, one of the following occurs: \\
\begin{enumerate}
\item An $l$-thin curve $v$ splits the pair $(\Sigma_i ,\Sigma_{i+1})$, i.e. $v$ splits the associated split level surfaces $(S_i ,S_{i+1})$, which in turn  form
an $l$-thin pair. 
\item there exists an $L_1$-bi-Lipschitz embedding  $$G_i: (S\times [0,1], (\partial S)\times [0,1]) \rightarrow (M, \partial M)$$
such that $\Sigma_i^s = G_i (S\times \{ 0\})$ and $\Sigma_{i+1}^s = G_i (S\times \{ 1\})$
\end{enumerate}
Finally, each $l$-thin curve in $S$ splits at most
$n$ split level surfaces in the  sequence $\{ \Sigma_{i} \}$. \label{wsplit}
\end{theorem}

A model manifold $M$ all of whose ends are equipped with a collection of  exiting    split  surfaces satisfying the conclusions of Theorem \ref{wsplit} is said to be equipped
with a {\bf weak split geometry} structure.

Pairs of split surfaces satisfying Alternative (1) of Theorem \ref{wsplit} will be called an $l$-thin pair of split surfaces (or simply a thin
pair if $l$ is understood). Similarly, pairs of split surfaces satisfying Alternative (2) of Theorem \ref{wsplit} will be called an $l$-thick pair
(or simply a thick
pair) of split surfaces.

\begin{defn} Let $(\Sigma_i^s, \Sigma_{i+1}^s)$ be a thick pair of split surfaces in $ M$. 
The closure of the bounded component of
$M \setminus (\Sigma_i^s \cup \Sigma_{i+1}^s)$ between  $\Sigma_i^s, \Sigma_{i+1}^s$   will be called a thick block.\end{defn}

Note that a thick block is uniformly bi-Lipschitz to the product $S \times [0,1]$ and that its boundary components are 
$\Sigma_i^s, \Sigma_{i+1}^s$.

\begin{defn} \label{splitb} Let $(\Sigma_i^s, \Sigma_{i+1}^s)$ be an $l$-thin pair of split surfaces in $M$
and $F(\TT_i)$ be the collection of $l$-thin Margulis tubes that split both $\Sigma_i^s, \Sigma_{i+1}^s$. The closure of the union of the
bounded components of
$M \setminus ((\Sigma_i^s \cup \Sigma_{i+1}^s)\bigcup_{F(T)\in F(\TT_i)} F(T))$  between  $\Sigma_i^s, \Sigma_{i+1}^s$  will be called a split block.
Equivalently, the closure of the union of the
bounded components of
$M(l) \setminus (\Sigma_i^s \cup \Sigma_{i+1}^s)$  between  $\Sigma_i^s, \Sigma_{i+1}^s$   is a split block. Each connected component of
 a split block is a split component.
\end{defn}

\begin{rmk}  \cite[Remark 4.12]{mahan-split} {\rm   For each lift $\til{K} \subset
\til{M}$ of a split component $K$ of a split block of $M(l) \subset M$, 
there are lifts of $l$-thin Margulis tubes that share the boundary of  $\til{K}$ in $\til{M}$. Adjoining these lifts to 
$\til{K}$ we obtain {\bf extended split components}. Let $\KK^\prime$ denote the collection of extended split components  in $\til{M}$.
Denote  the collection of split components in $\til{M(l)} \subset \til{M}$ by $\KK$.
Let $\til{M(l)}$ denote the lift of $M(l)$ 
to $\til M$. 
Then the inclusion of $\til{M(l)}$ into $\til{M}$ gives a quasi-isometry between $\EE (\til{M(l)}, \KK)$ and
$\EE (\til{M}, \KK^\prime)$ equipped with the respective electric metrics. This  follows from the last assertion of Theorem \ref{wsplit}.

 The electric metric on $\EE (\til{M}, \KK^\prime)$ is called the {\bf graph-metric}  \cite[Section 4.3]{mahan-split} and
is denoted by $d_G$. The electric space will be denoted as $(\til{M}, d_G)$.

The electric metric on $\EE (\til{M}, \KK \bigcup \TT_l)$ is quasi-isometric to the electric metric on $\EE (\til{M}, \KK^\prime)$, again by  the last assertion of Theorem \ref{wsplit}. 
The electric space will be denoted as $(\til{M}, d_G^1)$.}
\end{rmk}

\begin{definition} Let $Y \subset \til{N}$ and $X=F(Y)$.  $X \subset \til{M}$ is said to 
be $\Delta$-graph quasiconvex if for any hyperbolic geodesic $\mu$ joining $a, b \in Y$,
$F(\mu )$ lies inside $N_\Delta (X, d_G) \subset \EE (\til{M}, \KK^\prime)$. \end{definition}

For $X$ a split component in a manifold, define $CH(X) = F(CH(Y))$, where $CH(Y)$
  is the convex hull of $Y$ in $\til{N}$, provided the ends of $N$ have no cusps, i.e. $N=N^h$. Else define $CH(X)$ to be the image under
$F$ of $CH(Y)$ minus cusps.
Further, in order to ensure hyperbolicity of the universal cover, we partially electrocute the cusps of $M$ (cf. Theorem \ref{pel}).

Then $\Delta$-graph quasiconvexity of $X$
is equivalent to the condition that  $dia_G (CH(X))$  is bounded by $\Delta^{\prime} = \Delta^{\prime}(\Delta )$ as any split component 
has diameter one in $(\til{M}, d_G)$. We recall the following from \cite{mahan-split}.

\begin{lemma}  \cite[Lemma 4.16]{mahan-split} Let $E$  be a simply 
 degenerate end of a simply or doubly degenerate hyperbolic 3-manifold $N$ homotopy equivalent to a surface $S$
and equipped with
a weak split geometry model $M$.
For $K$ a split component contained in $E$, let $\tilde{K}$ be a lift to $\til N$.  Then there exists $C_0
= C_0(K)$ such that  the convex hull of $\tilde K$ minus cusps lies in a
$C_0$-neighborhood of $\tilde K$ in $\til N$. \label{hypqc-a}
\end{lemma}

\begin{prop}  \cite[Proposition 4.23]{mahan-split} For $K$ a split component,  $\til{K} $ is
  uniformly graph-quasiconvex in  $\til{M}$, i.e. there exists $\Delta^\prime$ such that $dia_G(CH(\til{K}))) \leq  \Delta^\prime$
for all incompressible split components $\til{K} $.
\label{gr-qc-free-a}
\end{prop}

We summarize the conclusions of the above propositions below.

\begin{defn} \label{minsky-split} A model manifold of weak split geometry is said to be of {\bf split geometry} if \\
\begin{enumerate}
\item  Each split component $\til{K}$ is 
quasiconvex (not necessarily uniformly) in the hyperbolic metric on $\til{N}$. \\
\item Equip $\til{M}$ with the {\it graph-metric} $d_G$ obtained by
electrocuting  (extended) split components $\til{K}$. Then the convex hull
$CH( \til{K})$ of any split component $\til{K}$ has uniformly bounded
diameter in the metric $d_G$. 
\end{enumerate}
\end{defn}

Hence by Lemma \ref{hypqc-a} and Proposition \ref{gr-qc-free-a} we have the following technical Theorem of \cite{mahan-split}.

\begin{theorem}  \cite{minsky-elc1, minsky-elc2}  \cite[Theorem 4.32]{mahan-split}
Any simply or doubly degenerate  hyperbolic 3-manifold homotopy equivalent to a surface
 is bi-Lipschitz homeomorphic to a Minsky model
and hence to a model of split geometry. \label{gqc}
\end{theorem}

\subsubsection{Ladders:} For details on the construction of ladders, see \cite[Section 5]{mahan-split}.
Note  that after welding the boundary components of $S^s_i$ together in a split block,  we obtain a bounded geometry surface $S_i$
in $M_{wel}$. Thus $S_i$ is the {\it connected}  bounded geometry surface obtained from  $S^s_i$ by equipping it with the quotient topology dictated by  welding. 
For convenience of notation, we redesignate this surface $S_i$. In the  welded model manifold $M_{wel}$, we thus obtain a sequence of bounded geometry
surfaces $\{ S_i \}$ exiting the end(s). The region between $S_i$ and $S_{i+1}$ is either a thick block or a split block.

From a geodesic $\lambda = \lambda_0 \subset \tilde{S} \times \{ 0 \} 
\subset \widetilde{B_0}$ we constructed in \cite{mahan-split} a `hyperbolic ladder'
${\LL}_\lambda (\subset \til{M_{wel}})$ such that $\lambda_i  =  {\LL}_\lambda \cap  \tilde{S_i} $ is an
electro-ambient quasigeodesic
in the (path) electric metric on $ \tilde{S_i}$ induced by the graph metric $d_G$ on $\tilde{M}$. $\lambda_{i+1}$ is constructed inductively from $\lambda_i$
(in \cite{mahan-split} or \cite{mahan-amalgeo}) by `flowing $\lambda_i$ up' in the block $\widetilde{B_i}$. More precisely,
$\widetilde{B_i}$ has a natural product structure and is bounded by $\widetilde{S_i}$ and $\widetilde{S_{i+1}}$. Given $\lambda_i$
joining $p_i, q_i \in \widetilde{S_{i}}$, there exist points $p_{i+1}, q_{i+1} \in \widetilde{S_{i+1}}$ lying vertically above $p_i, q_i$
respectively. $\lambda_{i+1}$ is the electro-ambient geodesic in $\widetilde{S_{i+1}}$ (equipped with the electric metric) joining $p_{i+1}, q_{i+1}$.
 
We also constructed a large-scale retract $\Pi_\lambda : \tilde{M} \rightarrow {\LL}_\lambda$ such that
the restriction $\pi_i$ of $\Pi_\lambda$ to $\tilde{S} \times \{ i \}$ is, roughly speaking,
a nearest-point retract of $\tilde{S} \times \{ i \}$ onto $\lambda_i$ in the (path) electric 
metric on $ \tilde{S_i}$.

We have the following basic theorem from \cite{mahan-split}

\begin{theorem}  \cite[Theorem 5.7]{mahan-split}
There exists $C > 0$ such that 
for any geodesic $\lambda = \lambda_0 \subset \widetilde{S} \times \{
0 \} \subset \widetilde{B_0}$, the retraction $\Pi_\lambda :
\widetilde{M} \rightarrow {\LL}_\lambda $ satisfies: \\

Then $d_{G}( \Pi_{\lambda } (x), \Pi_{\lambda } (y)) \leq C
d_{G}(x,y) + 
C$.
\label{retract}
\end{theorem}

\subsubsection{qi Rays:}
We also have the following from \cite{mahan-split}.

\begin{lemma} \cite[Lemma 5.9]{mahan-split}
There exists $C \geq 0$ such that for $x_{i} \in \lambda_{i}$ there
exists
$x_{i-1} \in \lambda_{i-1}$ with $d_G (x_{i}, x_{i-1}) \leq C$. Similarly 
there
exists
$x_{i+1} \in \lambda_{i+1}$ with $d_G (x_{i}, x_{i+1}) \leq C$. Hence, for all $n$ and $x \in \lambda_n$,
 there exists a $C$-quasigeodesic ray $r$ such that $r(i) \in \lambda_i \subset {\LL}_\lambda$ for all $i$
and $r(n) = x$. 
\label{dGqgeod}
\end{lemma}

Further, by construction of split blocks, $d_G (x_{i}, S_{i-1}) =
1$. Therefore inductively, $d_G (x_{i}, S_{j}) =
|i-j|$. Hence $d_G (x_{i}, x_{j}) \geq
|i-j|$. By construction, $d_G (x_{i}, x_{j}) \leq
C|i-j|$.

Hence, given $p \in \lambda_{i}$ the sequence of points $ x_{n}, n \in \mathbb{N} \cup \{ 0 \}$ (for simply degenerate groups)
or $n \in \integers$ (for totally degenerate groups) with $x_i = p$ gives by Lemma \ref{dGqgeod} above, a quasigeodesic in
the $d_G$-metric. Such quasigeodesics shall be referred to as {\em $d_G$-quasigeodesic rays}.

\section{Laminations}\label{lamns}
\subsection{Ideal points  are identified by  Cannon-Thurston Maps }\label{easy}
We would like to know exactly which points are identified by the
Cannon-Thurston map, whose existence is ensured by Theorem \ref{split}. Let $i: \til{S} \rightarrow \til{M}$ denote inclusion.
Let $\hat{i}$ be the continuous extension of $i$
to the disk ${\mathbb{D}} = ({\mathbb{H}}^2 \cup S^1_\infty)$ in Theorem \ref{split}.
Let $\partial i$ denote the restriction of $\hat{i}$ to the boundary $S^1_\infty$.

As mentioned in the introductory Section \ref{out}, we shall first prove the {\bf forward direction}
of Theorem  \ref{main}.
 Proposition \ref{ptpreimage} below shows that the existence of
a Cannon-Thurston map automatically guarantees  that end-points of
leaves of the ending lamination are identified by the Cannon-Thurston
map.

\begin{prop}
 Let $u, v$ be either ideal 
end-points of a leaf of an ending lamination, or ideal boundary points of a
complementary ideal polygon. Then
 $\partial i(u) = \partial i(v)$.
\label{ptpreimage}
\end{prop}

\begin{proof} ( cf. Lemma 3.5 of
 \cite{mitra-endlam}. See also \cite{mahan-split}.) We shall use two facts in the proof: \\
1) the fact (due to Bonahon \cite{bonahon-bouts} and  Thurston \cite{thurstonnotes}) that surface groups
are tame and that for $M$ there exist 
simple closed curves $a_i$ on $S$ whose geodesic realizations exit the end.\\
2) the fact (due to Thurston \cite{thurstonnotes} Ch. 9) that the sequence of simple closed curves $a_i$ converges to the ending lamination in the space of
measured laminations. It follows that after lifting to the universal cover,  any leaf of the ending lamination is a Chabauty topology  limit of bi-infinite
geodesics  $\til{a_i}$ (lifts of $a_i$).

Since any end $E$ of $M$ is geometrically tame \cite{thurstonnotes}, there exists $C_0$ such that
there exists a sequence of closed geodesics $s_i$ with length at most $C_0$
 exiting the
end. We shall refer to such geodesics as `bounded geodesics'. 
Let $a_i$ be geodesics in the intrinsic metric
on the base surface $S$ ( $=S_0 \subset M$)
freely homotopic to $s_i$.  We 
 can assume further \cite{bonahon-bouts} that $a_i$'s are simple. 
Join $a_i$ to
 $s_i$ by the shortest geodesic $t_i$
 in $M$ connecting the two
 curves. 

For any leaf $l$ of the ending lamination, we have a subsequence of the $a_i$'s whose Hausdorff limit in $S$ contains $l$.
Abusing notation slightly let us denote the subsequence as $\{ a_i \}$.
In the universal cover, we obtain segments
$a_{fi} \subset \til{S}$ 
which are finite segments whose end-points are identified by the
 covering map $P: \til{M} \rightarrow M$. We also assume that $P$ is
 injective restricted to the interior of $a_{fi}$'s mapping to
 $a_i$.
 Similarly there
 exist 
 segments
$s_{fi} \subset \til{M}$ 
which are finite segments whose end-points are identified by the
 covering map $P: \til{M} \rightarrow M$. We also assume that $P$ is
 injective restricted to the interior of $a_{fi}$'s. 
The finite segments $s_{fi}$ and $a_{fi}$ are chosen in such
 a way that there exist 
 lifts $t_{1i}$, $t_{2i}$, joining end-points of $a_{fi}$
 to corresponding end-points of $s_{fi}$. The union of these four pieces looks like a trapezium (See figure below, where we omit subscripts for convenience).

\begin{center}

\includegraphics[height=4cm]{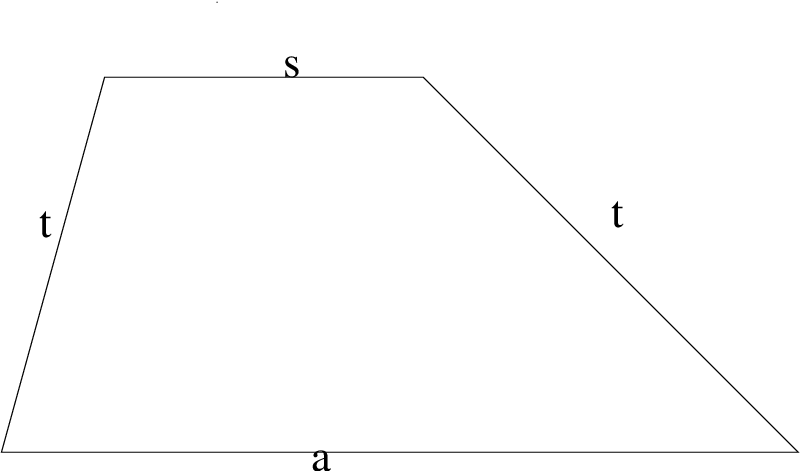}

\underline{Figure 2:  {\it Trapezium} }

\end{center}

\smallskip

Next, given any lift $\lambda$ of the leaf $l$ to $\til S$, we may choose
translates of the finite segments $a_{fi}$
(under the action of $\pi_1(S)$) appropriately, such that $a_{fi}$
converge to $\lambda$ in (the Hausdorff/Chabauty topology on closed subsets of) ${\Hyp}^2$. For each $a_{fi}$, let 

\begin{center}
$\beta_{fi} = t_{1i} \circ s_{fi} \circ \overline{t_{2i}}$
\end{center}
where $\overline{t_{2i}}$ denotes $t_{2i}$ with orientation
  reversed. Then $\beta_{fi}$'s are uniform hyperbolic quasigeodesics
  in $\til{M}$ (since $s_{fi} $ is short). If the translates of $a_{fi}$ we are considering
  have end-points lying outside large balls around a fixed reference
  point $p \in \til{S}$, it is easy to check that $\beta_{fi}$'s lie
  outside large balls about $p$ in $\til{M}$.

At this stage we invoke the existence theorem for Cannon-Thurston
maps, Theorem \ref{split}. Since $a_{fi}$'s converge to
$\lambda$ and there exist uniform hyperbolic quasigeodesics
$\beta_{fi}$, joining the end-points of $a_{fi}$ and exiting all
compact sets, it follows that  $\partial i(u) = \partial i(v)$, where
$a, b$ denote the boundary points of $\lambda$. 

Hence if we define $u, v$ to be equivalent if they are the end-points of a leaf of the ending lamination, then the transitive closure of this relation has as elements of an equivalence class \\
a) either ideal 
end-points of a leaf of a lamination, \\
b) or ideal boundary points of a
complementary ideal polygon, \\
c) or  a single point in $S^1_{\infty}$ which is not an end-point of a leaf of a lamination.
\end{proof}

\begin{definition} Let $H$ be a finitely presented group acting on a hyperbolic space $X$ with quotient $M$. Let $X_H$ be a 2-complex with fundamental group $H$, and $i:X_H \rightarrow M$ be a map inducing an isomorphism of fundamental groups. Then $i$ lifts to $\tilde{i}: \tilde{X_H} \rightarrow X$. A bi-infinite geodesic $\lambda$ in $\tilde{X_H} \subset X $ will be called a {\bf leaf of the abstract ending lamination} for $i: X_H \rightarrow M$, if \\
1) there exists a set of geodesics $\sigma_i$ in $M$ exiting every compact set \\
2) there exists a set of geodesics $\alpha_i$ in $X_H$ with $i(\alpha_i)$ freely homotopic to $\sigma_i$ \\
3) there exist finite lifts $\til{\alpha_i}$ of $\alpha_i$ in $(\til{X_H})$ such that the natural covering map $\Pi: \til{X_H} \rightarrow X_H$ is injective away from end-points of $\til{\alpha_i}$ \\
4) $\til{\alpha_i}$ converges to $\lambda$ in the Chabauty topology
\end{definition}

Proposition \ref{ptpreimage} and its proof readily generalize to

\begin{prop}
Suppose $H$ is hyperbolic and $\tilde{i}: \til{X_H} \rightarrow X$ extends to a Cannon-Thurston map on boundaries. 
 Let $u, v$ be 
end-points of a leaf of an abstract ending
lamination. Then
 $\partial i(u) = \partial i(v)$.
\label{ptpreimage1}
\end{prop}

To distinguish between the ending lamination and bi-infinite geodesics whose end-points are identified by $\partial i$, we make the following definition.

\begin{definition}
A {\bf CT leaf} $\lambda_{CT}$ is a bi-infinite geodesic whose end-points are identified by $\partial i$. \\ An {\bf EL leaf} $\lambda_{EL}$ is a bi-infinite geodesic whose end-points
are  ideal boundary points of
either a leaf of the ending lamination, or a complementary ideal polygon.
\end{definition}

Then to prove the main theorem \ref{main} it remains to show that  \\
$\bullet$ {\bf A {\it CT leaf} is an {\it EL leaf}.}

\subsection{Leaves of Laminations}\label{leaves}

Our first observation is that any semi-infinite geodesic (in the {\em hyperbolic metric} on $\tilde{S}$) contained in a {\it CT leaf} in the base surface $\tilde{S}$ $= \tilde{S} \times \{ 0 \} \subset \tilde{B} = \tilde{B}\times \{ 0 \}$ has infinite diameter in the {\em graph metric $d_G$} restricted to $\tilde{S} \times \{ 0 \}$, i.e. the
induced path metric on $\tilde{S} \times \{ 0 \}$. This follows from the following somewhat stronger assertion.

\begin{lemma} Given $k \geq 0$, there exists $C \geq 0$ such that if $B = \cup_{0\leq i \leq k} B_i$ and $\lambda \subset \til{B}$ is a bi-infinite
geodesic in the intrinsic metric on $\til{B}$, whose end-points are identified by the Cannon-Thurston map, then for any split component $\til{K}$, $dia_{hyp} (\lambda \cap 
\til{K} ) \leq C$
\label{findiacor}
\end{lemma}

\begin{proof} Suppose not. 

Then there exist split components $\til{K(i)} \subset \til{B}$, such that $dia_{hyp} (\lambda \cap 
\til{K(i)} ) \geq 3i$, where $dia_{hyp} $ denotes diameter in the hyperbolic metric on $\til{M}$. Acting on $\til{B}$ by elements $h_i$ of the surface group $\pi_1(S)$, we may assume that there exists a sequence of segments $\lambda^i \subset h_i \cdot \lambda$ such that \\
$\bullet $ $\lambda^i$ is approximately centered about a fixed origin $0$ in a fixed lift $\til{K}$ of a fixed split component $K$, i.e. $\lambda^i$ pass uniformly close to $0$ and end-points of $\lambda^i$ are at distance $\geq i$ from $0$. This is possible since $B$ contains finitely many split blocks. \\
Since $\til{K}$ is quasiconvex, it follows that the $\lambda^i$'s are uniform quasigeodesics in $\til{M}$. Hence, 
the sequence $\{ h_i \cdot \lambda \}$ converges to
 a bi-infinite quasigeodesic $\lambda^\infty$  in the Chabauty topology. Since the set of {\em CT leaves}
are closed in the Chabauty topology, it follows that $\lambda^\infty$ is a {\em CT leaf}. 

But, this is a contradiction, as  we have noted already that $\lambda^\infty$ is a quasigeodesic. 
\end{proof}

\begin{cor} {\bf CT leaves have infinite diameter} \\
Let $\lambda_+ ( \subset \lambda \subset \tilde{S} \times \{ 0 \} = \tilde{S})$ be a semi-infinite geodesic (in the {\em hyperbolic metric} on $\tilde{S}$) contained in a {\it CT leaf} $\lambda$. Then $dia_G (\lambda_+)$ is infinite, where $dia_G$ denotes diameter in the graph metric restricted to 
$\til{S}$.
\label{ctleafinf}
\end{cor}

\begin{proof}  Put $k=1$ in Lemma \ref{findiacor}. \end{proof}

Using Lemma \ref{findiacor}, we shall now show:

\begin{prop}
There exists a function $M(N) \rightarrow \infty$ as $N \rightarrow \infty$ such that the following holds: \\
Let $\lambda$ be a CT leaf. Also for $p, q \in \til{M}$, let $\overline{pq}$ denote a geodesic in $(\til M, d_G)$ joining $p, q$.
If $a_i, b_i \in \lambda$ be such that $d(a_i,0) \geq N$, $d(b_i,0) \geq N$, then $d_G ({\overline{a_i b_i}}, 0) \geq M(N)$, where $d$ denotes the hyperbolic metric on $\til{M}$ and $d_G$ the graph metric.
\label{ctdg}
\end{prop}

\begin{proof} Suppose not. Let $\lambda_+$ and $\lambda_-$ denote the ideal end points of $\lambda$. 
 Then there exists $C \geq 0$, $a_i \rightarrow \lambda_-$, $b_i \rightarrow \lambda_+$ such that $d_G ({\overline{a_i b_i}}, 0) \leq C$. That is, there exist $p_i \in {\overline{a_i b_i}}$ such that $d_G(0,p_i) \leq C$. Due to the existence of a Cannon-Thurston map in the hyperbolic metric (Theorem \ref{split}), we may assume that $d(0,p_i) \geq i$ (in the hyperbolic metric). Then the hyperbolic geodesic $\overline{0,p_i}$ passes through at most $C$ split blocks 
(cf. Definition \ref{splitb}) for every $i$. Let $B = \cup_{0\leq i \leq C} B_i$ and $p_i \rightarrow p_\infty$. Then $\overline{0,p_i} \subset \til{B}$. But since
$p_i \in {\overline{a_i b_i}}$,  the Cannon-Thurston map identifies 
$\lambda_-, \lambda_+, p_\infty$. See Figure below.

\begin{center}

\includegraphics[height=3.5cm]{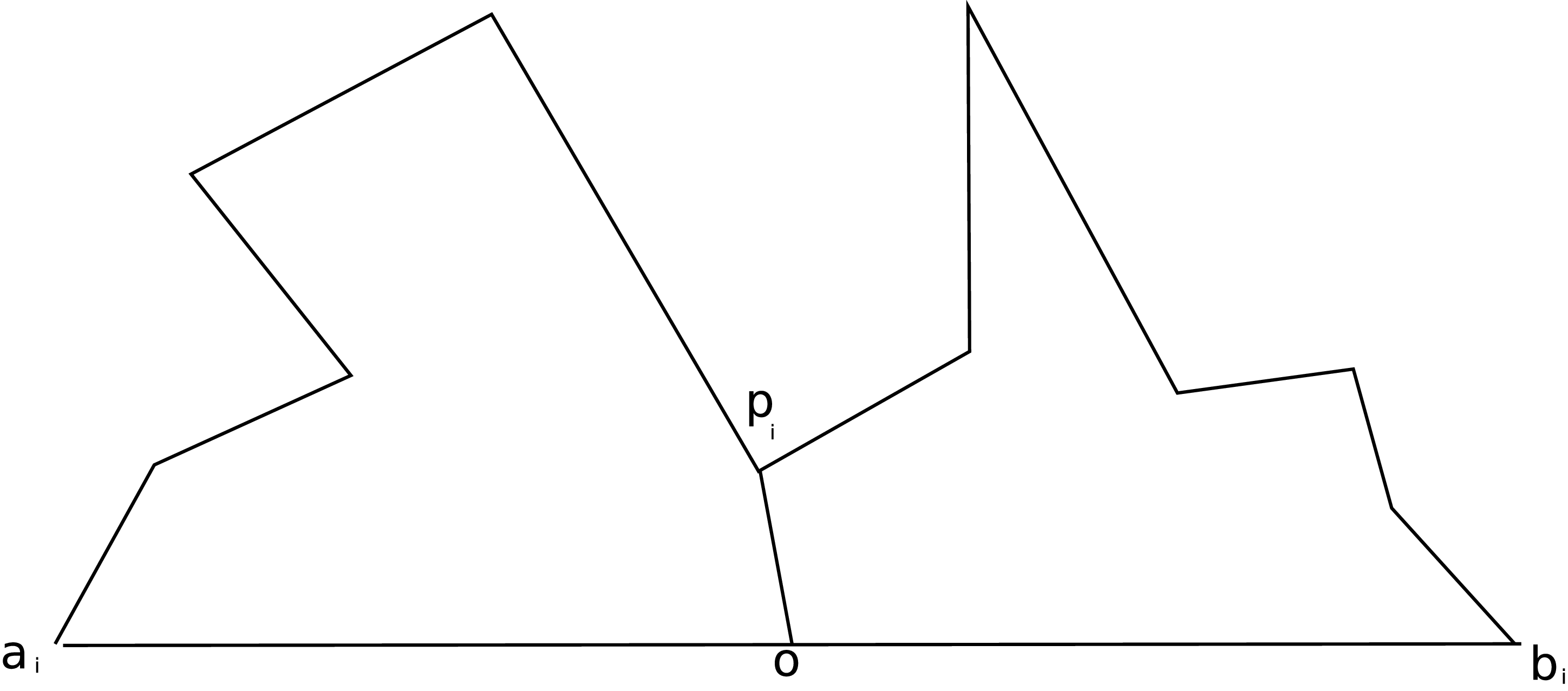}

\underline{Figure 3:  {\it Cannon-Thurston in the Graph Metric} }

\end{center}

\smallskip
Also, \\
\begin{center}
$\overline{0,p_\infty} \subset \overline{0,p_\infty} \cup \overline{0,\lambda_+}$\\
$\overline{0,p_\infty} \subset \overline{0,p_\infty} \cup \overline{0,\lambda_-}$\\
\end{center}

\noindent and at least one of the above two ($ \overline{0,p_\infty} \cup \overline{0,\lambda_+} = \overline{p_\infty,\lambda_+}$ say)
must pass close to $0$. Then $\overline{p_\infty,\lambda_+}$ is a {\em CT leaf}. But $ \overline{0,p_\infty} $ lies in a $C$-neighborhood of $0$ in the graph-metric $d_G$, contradicting Lemma \ref{findiacor} above.
\end{proof}

\section{Closed Surfaces}\label{closed}
In this section $S$ will denote a closed surface. As mentioned in the introductory Section \ref{out}, we shall now proceed to
 prove the {\bf reverse direction}
of Theorem  \ref{main}. The aim of this Section is to show that a {\bf CT} leaf is an {\bf EL} leaf.

\subsection{Geodesic Laminations and $\mathbb{R}$-trees} \label{S:laminations}

For a discussion of geodesic laminations (or simply
laminations as we shall call them), 
we refer the reader to \cite{harer-penner}, \cite{CEG}, 
\cite{thurstonnotes}, \cite{CB}. For a discussion on dual 
$\mathbb{R}$-trees, see \cite{shalen-dend}.

The space of filling laminations which we denote $\FL$
 are the measure classes of
measured laminations $\Lambda$ for which all complementary regions of 
the support $|\Lambda |$ are simply connected. The quotient of $\FL$ by
forgetting the measures will be denoted $\EL$ and is the space of 
\textit{ending laminations}. It is a well-known fact \cite{thurstonnotes, minsky-elc1} that ending laminations have no
simple closed leaves. 
A useful fact is that such
 laminations  are
minimal, i.e. the closure
(in the Hausdorff topology) of {\bf any} of its leaves
is the whole lamination.  We can identify a minimal lamination $\Lambda$ with
a closed invariant (under $\pi_1(S)$)
subset of the set of {\it unordered pairs} in 
$(S^1_\infty \times S^1_\infty \setminus \Delta ) /{\bf R_a}$, where
$\Delta$ denotes the diagonal and ${\bf R_a}$
is the relation identifying $(a,b)$ with $(b,a)$.

\begin{lemma}
Let $\Lambda$ be a minimal geodesic lamination on a surface $S$. Let $I$ be an embedded (closed)
interval in $S$ transverse to $\Lambda$. Let $\til{\Lambda}$ denote the
union of all lifts of leaves of $\Lambda$ to the universal
cover $\til{S}$. Let $\til{I}$ denote the union of all lifts
of $I$ to $\til{S}$. Define two leaves of $\til{\Lambda}$ to be equivalent
if both of them intersect the same component of $\til{I}$.
Then the limit set of any connected component of
the transitive closure of this relation contains a pair of poles $g^\infty$ and $g^{-\infty}$ for some element $g \in \pi_1(S)$.
\label{minimal}
\end{lemma}

\begin{proof} Let $\TT$ be the $\mathbb{R}$-tree dual to
$\til{\Lambda}$. Let $I_0$ be a fixed lift of $I$ to $\til{S}$. Then $I_0 \subset \til{S}$ 
projects to an embedded non-trivial interval (also called $I_0$)
in $\TT$ under the quotient map $q$ that identifies leaves
of $\til{\Lambda}$ to points. 
The orbit of $I_0$ under $\pi_1(S)$ acting on $\TT$
is a forest, in fact a sub-forest $\FF$ of $\TT$.

Let  $\TT_1$ be  the connected component of $\FF$ containing $I_0$. 
If $g \TT_1 \cap \TT_1 \neq \emptyset$, then $g \TT_1 \subset \TT_1$.
Hence $g^{-1}\TT_1 \cap \TT_1 \neq \emptyset$, and we finally 
have that $\TT_1$ is invariant under $g^n$ for all integers $n$.
This shows that $q^{-1} (\TT_1) \subset \til{S}$ contains the
pole corresponding to the infinite order element $g$.

Thus we need finally the existence of a $g$ as in the previous paragraph.
It suffices to show that for any non-trivial $I_0$, there exists
$g \in \pi_1(S)$ such that $gI_0\cap I_0 \neq \emptyset$. But this 
follows from minimality of $\Lambda$, using the fact that
each leaf is dense in $\Lambda$, and hence that there exists
$g \in \pi_1(S)$ such that $gI_0$ and $I_0$ are transverse to
a {\it common leaf} $\lambda$ of $\til{\Lambda}$. 
\end{proof}

\subsection{Rays Contained in Ladders}

\begin{definition}  Let $X, Y, Z$ be geodesically complete metric
spaces such that $X \subset Y \subset Z$. $X$ is said to {\bf coarsely separate}
$Y$ into $Y_1$ and $Y_2$ if \\
(1)  ${Y_1}\cup{Y_2} = Y$ \\
(2)  ${Y_1}\cap{Y_2} = X$ \\
(3)  For all $M \geq 0$, there exist ${y_1}\in{Y_1}$ and 
${y_2}\in{Y_2}$ such that $d({y_1},{Y_2}) \geq M$ and 
$d({y_2},{Y_1}) \geq M$  \\
(4)  There exists $C \geq 0$ such that for all ${y_1}\in{Y_1}$ 
and ${y_2}\in{Y_2}$ any geodesic in $Z$ joining ${y_1}, y_2$ passes through
a $C$-neighborhood of $X$. 
\label{cs}
\end{definition}

Let $\lambda = \lambda_0$ be any bi-infinite geodesic in $\til{S}$. Let $\LL_\lambda$ be the ladder corresponding to $\lambda$ as in Theorem  \ref{retract}. 

We now fix a quasigeodesic ray $r_0$ as in Lemma \ref{dGqgeod}, and consider a translate $r^{\prime} = h\cdot r_0$ passing through $z \in \lambda_m \subset {\LL}_\lambda$, i.e. $ r^{\prime}(m) = z$. Let $\Pi_\lambda \cdot 
r^{\prime} = r \subset {\LL}_\lambda$. Each $r(i)$ cuts $\lambda_i$ into two pieces $\lambda_i^-$ and $\lambda_i^+$ with ideal boundary points 
$\lambda_{i,-\infty }, \lambda_{i,\infty }$ respectively. 

We shall show that $r$ coarsely separates ${\LL}_\lambda$ into ${\LL}_\lambda^{+ }$ and ${\LL}_\lambda^{-}$, where 
\begin{center}
${\LL}_\lambda^{+ } = \bigcup_i \lambda_i^+$ \\
${\LL}_\lambda^{- } = \bigcup_i \lambda_i^-$ \\
\end{center}

\noindent and $\lambda_i^+$ (resp. $\lambda_i^-$ ) is the segment of $\lambda$ joining $r(i)$ to the ideal end-point $\lambda_{i,-\infty }$ (resp. $\lambda_{i,\infty }$. 

We need to repeatedly apply Theorem \ref{retract} to prove the above assertion. 

Given $r^{\prime}$, we construct two hyperbolic ladders 
${\LL}_\lambda^{+\prime }$ and ${\LL}_\lambda^{-\prime }$, obtained by joining the points $r^{\prime} (i)$ to the ideal end-points
 $\lambda_{i,-\infty }, \lambda_{i,\infty }$ respectively, of $\lambda_{i} \subset \til{S} \times \{ i \}$. Then 
${\LL}_\lambda^{+\prime }$ and ${\LL}_\lambda^{-\prime }$ are $C$-quasiconvex (in the graph metric $d_G$)
by Theorem \ref{retract}. Further, $\Pi_\lambda \cdot r^{\prime}(i) = r(i)$ by definition of $r$. Hence, \\

\begin{center}
$\Pi_\lambda ({\LL}_\lambda^{-\prime }) = {\LL}_\lambda^{-} $\\
$\Pi_\lambda ({\LL}_\lambda^{+\prime }) = {\LL}_\lambda^{+}$ \\
\end{center}

Further, \\
\begin{center}
${\LL}_\lambda^{-}  \cup  {\LL}_\lambda^{+} = {\LL}_\lambda  $\\
${\LL}_\lambda^{-} \cap {\LL}_\lambda^{+} = r $ \\
\end{center}

\noindent and there exists a $K_0$ (independent of $r_0, h, \lambda$) such that 
${\LL}_\lambda^{-} , {\LL}_\lambda^{+} , {\LL}_\lambda, r$ are all $K_0$-quasiconvex. 

Criterion (3) of Definition \ref{cs} in this context is given by Lemma \ref{ctleafinf} :  {\it CT leaves have infinite diameter}.

To prove that $r$ separates $ {\LL}_\lambda$ into ${\LL}_\lambda^{-} , {\LL}_\lambda^{+} $, we need to show first:

\begin{lemma} 
For all $K_0 \geq 0$, there exists $K_1 \geq 0$ such that if $p \in 
{\LL}_\lambda^{-} , q \in {\LL}_\lambda^{+} $ with $d_G(p,q) \leq K_0$, then there exists $z \in r$ such that $d_G(p,z) \leq K_1$ and $d_G(q,z) \leq K_1$.
\label{cs2}
\end{lemma}

\begin{proof} Let $\Pi_\lambda^+$ denote the sheetwise retract of Theorem \ref{retract} onto ${\LL}_\lambda^+$. Then 
$$\Pi_\lambda^+ (\lambda_i^- ) = r(i)$$ 

and 
  $$\Pi_\lambda^+ (x) = x$$ 
for all $x \in {\LL}_\lambda^+$.

Hence

\begin{center}
$\Pi_\lambda^+ (q) = q $\\
$\Pi_\lambda^+ (p) =  z = r(i)$
\end{center}
for some $z \in r$ and some $i$. 

Therefore, by Theorem \ref{retract} again,
$$d_G(q,z) \leq Cd_G(p,q) = CK_0.$$

\noindent Choosing $K_1 = CK_0 + K_0$ (and using the triangle inequality for $p,q,z$) the Lemma follows.
\end{proof}

We are now in a position to prove:

\begin{theorem}
$r$ coarsely separates $ {\LL}_\lambda$ into ${\LL}_\lambda^{-} , {\LL}_\lambda^{+} $.
\label{csthm}
\end{theorem}

\begin{proof} We have already shown \\
\begin{center}
${\LL}_\lambda^{-} \cup {\LL}_\lambda^{+} = {\LL}_\lambda$ \\
${\LL}_\lambda^{-} \cap {\LL}_\lambda^{+} = r$ 
\end{center}

\noindent and there exists a $K_0$ (independent of $r_0, h, \lambda$) such that 
${\LL}_\lambda^{-} , {\LL}_\lambda^{+} , {\LL}_\lambda, r$ are all $K$-quasiconvex. 

Criterion (3) of Definition \ref{cs}  is given by Lemma \ref{ctleafinf}. 

Finally given $u \in {\LL}_\lambda^{-}$ and $v \in {\LL}_\lambda^{+}$, let $\overline{uv}$ be the geodesic in $(\til{M}, d_G)$ joining $u, v$. Then $\Pi_\lambda
( \overline{uv} )$ is a "dotted quasigeodesic", i.e. there is a sequence of points
 $u=p_0, p_1, \cdots p_n = v$, where $d_G(p_i,p_{i+1}) \leq C$ (and the constant $C$ is obtained from Theorem \ref{retract}). 
Further, $p_0 \in {\LL}_\lambda^{-}, p_n \in {\LL}_\lambda^{+}$ and 
$p_i \in {\LL}_\lambda$ for all $i$. Therefore there exists $m$ such that 
$p_m \in {\LL}_\lambda^{-}, p_{m+1} \in {\LL}_\lambda^{+}$, with 
$d_G(p_m,p_{m+1}) \leq C$. Hence, by Lemma \ref{cs2},
there exists $K_1 \geq 0$ such that  there exists $z \in r$ with $d_G(p_m,z) \leq K_1$ and $d_G(p_{m+1},z) \leq K_1$.

Finally, by Theorem
 \ref{gqc}, $(\til{M}, d_G)$ is hyperbolic, and therefore the "dotted quasigeodesic" $u=p_0, p_1, \cdots p_n = v$ lies in a uniformly bounded neighborhood of the
geodesic $\overline{uv}$. That is,
there exists $C_1 \geq 0$ such that for all 
$u \in {\LL}_\lambda^{-}$ and $v \in {\LL}_\lambda^{+}$, the geodesic
$\overline{uv}$  in $(\til{M},d_G)$ joining $u, v$ passes through
a $C_1$-neighborhood of $r$. This proves (4) in Definition \ref{cs} and hence we conclude 
that $r$ coarsely $ {\LL}_\lambda$ into ${\LL}_\lambda^{-} , {\LL}_\lambda^{+} $. 
\end{proof}

We shall have need for the following Proposition, whose proof
is exactly along the lines of Theorem \ref{csthm} above.

\begin{prop}
Let $\mu$, $\lambda$ be two bi-infinite geodesics on $\til{S}$
such that $\mu \cap \lambda \neq \emptyset$. Then ${\LL}_\lambda
\cap \LL_\mu$ contains a quasigeodesic ray $r$ coarsely separating
both ${\LL}_\lambda$ and $\LL_\mu$.
\label{twin}
\end{prop}

\begin{remark} {\rm Proposition \ref{twin} generalizes readily to cusped surfaces $S^h$ to show that if 
$\mu$, $\lambda$ are two bi-infinite geodesics on $\til{S^h}$
such that $\mu \cap \lambda \neq \emptyset$, then ${\LL}_\lambda
\cap \LL_\mu$ contains a quasigeodesic ray $r$.\\
To see this it suffices to note that 
if $\mu \cap \lambda \neq \emptyset$, then $\mu_i \cap  \lambda_i \neq \emptyset$ for all $i$. Hence we may construct a 
quasigeodesic ray $r$ contained in both ${\LL}_\lambda$ and $\LL_\mu$. } \label{twin-cusp} \end{remark}

One last Proposition to be used in the proof of Theorem \ref{main} is the following which says in particular 
that any two quasigeodesic rays lying on ${\LL}_\lambda$ are asymptotic with respect to the graph metric $d_G$..

\begin{prop}{\bf Asymptotic Quasigeodesic Rays}\\
Given $K \geq 1$ there exists $\alpha$ such that 
if $\lambda$ is a CT-leaf   then there exists
$z\in\partial \til{M}$ satisfying the following:\\
If  $r_1$ and 
${r}_2$ are $K$-quasi-geodesic rays
contained in $\LL_{\lambda}$ then there exists
 $N \in \natls$ such that\\
1) $r_j(n) \rightarrow z = \partial i(\lambda_{-\infty }) =
\partial i(\lambda_{\infty })$ as $n \rightarrow \infty$, for $j=1,2$. \\
2) $d_G(r_1(n),r_2(n)) \leq \alpha$ for all $n \geq N$
\label{qgeodasymp}
\end{prop}

\begin{proof} By Proposition \ref{ctdg}, we find that if $a_i, b_i \in \lambda = \lambda_0$ such that $a_i, b_i$ converge to ideal points $\lambda_{0,-\infty }, \lambda_{0,\infty }$ (denoted $\lambda_{-\infty }, \lambda_{\infty }$ for convenience), then $\Pi_\lambda (\overline{a_ib_i})$ leaves large balls about $0$. 
More precisely there exists $L_i \rightarrow \infty$ as $i \rightarrow \infty$ such that  $\Pi_\lambda (\overline{a_ib_i})$ lies outside the $L_i-$ball about $0$.

Also, by Theorem \ref{csthm} above, each $r_j$ coarsely
separates ${\LL}_\lambda$.
 Hence $\Pi_\lambda (\overline{a_ib_i})$ passes close to $r_j(n_{i(j)})$ for some 
$n_{i(j)} \in \natls$, where $n_i \rightarrow \infty$ as $i \rightarrow \infty$. We conclude that any such $r_j$ converges on $\partial \til{M}$ to the same point as 
$\partial (\lambda_{-\infty } ) = \partial ( \lambda_{\infty })$. This proves (1).

In particular any two quasigeodesic rays lying on ${\LL}_\lambda$ are asymptotic with respect to the graph metric $d_G$. This proves (2). 
\end{proof}

\subsection{Main Theorem for Simply Degenerate Groups} \label{sec:sd}
We are now in a position to prove the main Theorem \ref{main} of this paper for closed surfaces. For ease of exposition we shall deal with the simply degenerate case first
and then indicate the additional niceties for doubly degenerate groups. Recall that for a simply degenerate manifold $M = S \times J$, where $J = [0, \infty )$. For a totally
degenerate manifold $J = (-\infty , \infty )$ and it is the presence of two ends, positive and negative, that necessitates further care. The split level 
surfaces are indexed by $0, 1, \cdots , \infty$ for a simply degenerate manifold and by $\integers$ for  a totally (doubly) degenerate manifold. For a doubly degenerate
group, there will be {\it two} ending laminations, one for each end and a slight modification of the proof below will be necessary to identify and distinguish these. 
The constructions of ladders and blocks are otherwise identical in both cases.

\smallskip

\begin{theorem} \label{sdmain} 
 Let $\partial i(a) = \partial i(b)$ for $a, b \in
S^1_\infty$ be two distinct points that are identified by the
Cannon-Thurston map corresponding to a simply degenerate closed surface
group (without accidental parabolics). 
Then $a, b$ are either ideal 
end-points of a leaf of the ending lamination (in the sense of Thurston), 
or ideal boundary points of a
complementary ideal polygon. Further, if $a, b$ are either ideal 
end-points of a leaf of a lamination, or ideal boundary points of a
complementary ideal polygon, then $\partial i(a) = \partial i(b)$. \end{theorem}

\medskip

\begin{proof} The second statement has been shown in Proposition \ref{ptpreimage}. \\
To prove the first  statement, let $\partial i(a) = \partial i(b)$ for $a, b \in
S^1_\infty$. Then $(a,b) = \lambda \subset \til{S_0} \subset \til{M}$ 
is a CT-leaf.

Suppose $\lambda$ and $\mu$ are intersecting CT leaves, i.e.
$\partial i(\lambda_{-\infty }) =
\partial i(\lambda_{\infty })$  and $\partial i(\mu_{-\infty }) =
\partial i(\mu_{\infty })$. 

 As before, let
$\lambda_i$ and $\mu_i$ be
intersections of the ladders ${\LL}_\lambda$ and ${\LL}_\mu$
with the horizontal sheets. Then $r(i) = \lambda_i \cap \mu_i$ is a 
quasigeodesic ray by Proposition \ref{twin}. 
By Proposition \ref{qgeodasymp}, $r(i)$ converges
to a point $z$ on $\partial \til{M}$ as $i \rightarrow \infty$
such that $z =  \partial i(\lambda_{-\infty }) =
\partial i(\lambda_{\infty })  = \partial i(\mu_{-\infty }) =
\partial i(\mu_{\infty })$. Hence  the Cannon-Thurston map identifies the endpoints of any
two intersecting CT leaves $\lambda$
 and $\mu$.

If possible, suppose that the CT leaf $\lambda$ is \underline{not}
an EL-leaf. Then $\lambda$ intersects the ending lamination
transversely (since the ending lamination is a filling lamination
without any closed leaves) and there exist EL-leaves 
$\mu$ for which $\lambda \cap \mu \neq \emptyset$.
By Proposition \ref{ptpreimage}, each such $\mu$
is a CT-leaf. Hence, by the previous paragraph,
  the Cannon-Thurston map $\partial i$ identifies the end points of $\lambda$
with the endpoints of each such EL-leaf
$\mu$ for which $\lambda \cap \mu \neq \emptyset$. Let $z (\in S^2_\infty)$ denote this common image under  $\partial i$.

Since $\lambda$ is not  an EL-leaf, it contains a non-trivial
geodesic subsegment $I$ transverse to the ending lamination. Then the common image (under  $\partial i$) of end-points of all EL leaves $\mu$ intersecting $I$
transversely is $z$.

\begin{comment}
Let $\Lambda$  denote the
union of all lifts of leaves of the ending lamination to the universal
cover $\til{S}$. Let $\til{I}$ denote the union of all lifts
of $I$ to $\til{S}$. Define two leaves of $\til{\Lambda}$ to be equivalent
if both of them intersect a single element of $\til{I}$.
Then by Lemma \ref{minimal}
the limit set $\mathbb{L}$ of any equivalence class of
the transitive closure of this relation
\end{comment} 
By Lemma  \ref{minimal}  $(\partial i)^{-1}(z)$  contains a pair of poles $g^{-\infty} , g^\infty$ for some $g \in \pi_1(S)$. This is because the
equivalence class defined by $I$ as in Lemma  \ref{minimal} consists of pairs of points all of which are identified  (under  $\partial i$)
with $z$. 

\begin{comment}
All the points of $\mathbb{L}$ are identified to a point by the existence
of Cannon-Thurston maps
(Theorem \ref{split}) and by Proposition \ref{ptpreimage}.
Hence, in particular the Cannon-Thurston map identifies a pair of poles
to a point. 
\end{comment} 
This is a contradiction as a pair of poles forms the end-points
of a quasigeodesic in $\til{M}$.
We conclude that $\lambda$ must be
an EL-leaf. 
\end{proof}

\subsection{Modifications for Totally Degenerate Groups} We elaborate on  the 
modifications  indicated in the first paragraph of Section \ref{sec:sd}, to pass from the simply degenerate case to the totally degenerate case.
The construction of the `hyperbolic ladder' ${\LL}_\lambda$ as in the discussion preceding Theorem \ref{retract}
is done with indexing set $\mathbb{Z}$ in place of $\mathbb{N}$. In particular the quasigeodesic ray of
Lemma \ref{dGqgeod} is replaced by a bi-infinite quasigeodesic $r$. However, as a hyperbolic
metric space $\mathbb{Z}$ has two boundary points $+ \infty$ and $-\infty$. Correspondingly we have
two ending laminations $\Lambda_+$ and $\Lambda_-$. The easy direction of Theorem \ref{main}
given by Proposition \ref{ptpreimage} then goes through verbatim to show that
$\Lambda_+ \cup \Lambda_- \subset \Lambda_{CT}$. 

We need to find a way of distinguishing the
$+$ and $-$ directions in $\Lambda_{CT}$. To implement this, note that the discussion preceding
Proposition
\ref{qgeodasymp} shows that if $\lambda \in \Lambda_{CT}$, i.e. $\partial i (\lambda_\infty )
= \partial i (\lambda_{-\infty })$, then we have a bi-infinite quasigeodesic
$r: \mathbb{Z} \rightarrow {\LL}_\lambda$ such that $\partial i (\lambda_\infty )
= \partial i (\lambda_{-\infty }) = r(\alpha )$, where $\alpha$ is either $+\infty$ or $-\infty$.
Define $\Lambda_{CT}^+\subset\Lambda_{CT}$ (resp. $\Lambda_{CT}^- \subset\Lambda_{CT}$) to be
the collection of $CT$-leaves whose endpoints are identified in the $+\infty$ (resp. $-\infty$)
direction, i.e. $\alpha =$ $+\infty$ (resp. $-\infty$). Then 
the forward direction of Theorem \ref{main}
given by Proposition \ref{ptpreimage}  shows that
$\Lambda_+  \subset \Lambda_{CT}^+$ and $\Lambda_-  \subset \Lambda_{CT}^-$. Since both ending laminations
$\Lambda_+$ and $\Lambda_-$ are individually filling arational minimal laminations, the proof of Theorem
\ref{sdmain} (the reverse direction for simply degenerate groups) now shows that in fact 
$\Lambda_+  = \Lambda_{CT}^+$ and $\Lambda_-  = \Lambda_{CT}^-$.

\subsection{Application: Rigidity} In \cite{minsky-elc2}, Brock-Canary-Minsky prove the following Rigidity Theorem.

\begin{theorem} \label{bcmrig}
Let $G$ be a  closed surface group.
 If $\rho$ and $\rho'$ are two discrete faithful representations of
  $G$ into $PSl_2(\mathbb{C})$ that are conjugate by an
  orientation-preserving homeomorphism of $\hhat{\mathbb{C}}$, then $\rho$ and $\rho'$ are
  quasiconformally conjugate. \end{theorem}

We strengthen this by weakening the hypothesis of Theorem \ref{bcmrig} to a topological conjugacy only on limit sets (rather than all of $\hhat{\mathbb{C}}$).

\begin{theorem} \label{rig}
Let $G$ be a  closed surface group.
Let $\rho (G) = \Gamma$ and $\rho_1 (G) = \Gamma_1$ be two  simply or doubly degenerate  representations of
  $G$ into $PSl_2(\mathbb{C})$ with limits sets $\Lambda, \Lambda_1$. Suppose that the $G-$ actions on $\Lambda, \Lambda_1$ are topologically conjugate.
Then   $\rho$ and $\rho_1$ are
  quasiconformally conjugate. \end{theorem}

\begin{proof} 
We first deal with the simply degenerate case.
By Theorem \ref{main}, the pre-images of the Cannon-Thurston maps $\partial i$ and $\partial i_1$
from $\partial G (= S^1)$ to $\Lambda$ or $\Lambda_1$ are given by
end-points of leaves of the ending lamination (or ideal points of complementary polygons) whenever $\partial i$ and $\partial i_1$ are non-injective.
Thus the $G-$ action on $\Lambda, \Lambda_1$ pulls back to a $G-$ equivariant homeomorphism $\phi : \partial G \rightarrow \partial G$ taking 
the ending lamination of $\rho$ to that of $\rho_1$. Re-marking by an isomorphism of $G$ if necessary, we can ensure that the homeomorphism 
 be the identity on $\partial G$.
Hence the ending laminations of $\rho$ and $\rho_1$ are the same. 

In the doubly degenerate case, the same argument shows that the pairs of ending laminations for $\rho$ and $\rho_1$ are the same. 
Hence if $\rho$ and $\rho_1$ are doubly degenerate, they have the same end-invariants. By the Ending Lamination Theorem \cite{minsky-elc2},
$\rho$ and $\rho_1$ are {\bf conformally conjugate}.

In the simply degenerate case, the conformal structures corresponding to the geometrically finite ends for
$\rho$ and $\rho_1$  are quasiconformal deformations of each other (since the quotient of the domain of discontinuity is a connected
finite volume Riemann surface).
Since the ending laminations of $\rho$ and $\rho_1$  are the same, it follows therefore that the quotient manifolds are bi-Lipschitz homeomorphic
by the Ending Lamination Theorem \cite{minsky-elc2}. Hence $\rho$ and $\rho_1$ are
{\bf  quasiconformally conjugate}.
\end{proof}

\bigskip

\appendix

\Appendix{(by Shubhabrata Das and Mahan Mj) Surfaces with Cusps}\footnote{The work in this Appendix forms part of SD's PhD thesis written under the supervision of MM.
The proof given here
was discovered jointly considerably after the work on the earlier Sections was completed. Hence we have retained both approaches.}
\label{appendix}

We now deal with surfaces with cusps. $S^h$ will denote a finite volume hyperbolic surface with cusps. $S$ will denote
a truncated surface, i.e. $S^h$ minus an open neighborhood of the cusps. 
The arguments in this Section can be easily adapted to furnish a slightly different proof of Theorem \ref{main} for surfaces
without cusps. 
As in the previous Section we deal with the case of simply degenerate groups first.

\medskip

\noindent {\bf  Equivalence Relations on $S^1$:} \\ Suppose that a group $G$ acts on $S^1$ preserving a closed
 equivalence relation $\mathcal R$. An example $\RR_\LL$ of such a relation comes from a lamination $\LL$, where two points on $S^1$
are declared equivalent if they are end-points of a leaf  of $\LL$ .
 The equivalence relation $\RR_\LL$ is obtained as the transitive closure of this relation.

\begin{defn} \cite{bowditch-ct} Two disjoint subsets, $P, Q \subset S^1$ are linked if there exist linked pairs, $\{x, y \} \subset P$ and $\{ z, w \} \subset Q$.
$\RR$ is unlinked if distinct equivalence classes are unlinked. \end{defn}

The  following Lemma due to Bowditch give us a way of recognizing relations coming from laminations.

\begin{lemma} (Lemma 9.2 of \cite{bowditch-ct}) Let $\RR$ be a non-empty closed unlinked $G$-invariant equivalence relation on
$S^1$. Suppose that no pair of fixed points of any loxodromic are identified under $\RR$. Then
there is a unique complete perfect lamination, $\LL$, on $S$ such that $\RR = \RR_\LL$. \label{unlink} \end{lemma}

Let $\RCT$ denote the equivalence relation on $S^1$ induced by the Cannon-Thurston map for a simply degenerate
 punctured surface group (cf. Theorem \ref{split}).
Let $\Lambda$ denote the ending lamination. 
By Proposition \ref{ptpreimage} pairs of end-points of leaves of $\Lambda$ are contained in $\RCT$. Hence, for simply degenerate groups,
 it suffices to show that 
$\RCT$ is induced by a lamination since no other lamination can properly contain $\Lambda$. By Lemma \ref{unlink}
it suffices to show that $\RCT$ is unlinked. The next Proposition is the analogue of Proposition \ref{qgeodasymp} for cusped surfaces.

\begin{prop} Let $i: {\til {S^h}} \rightarrow \til{M^h}$ be an inclusion of the universal cover of a punctured surface
into the universal cover of the convex core $M^h$ of a simply degenerate 3-manifold. Let $\partial i$ be the associated Cannon-Thurston map.
If $\lambda$ is a CT-leaf in $\til S^h$, $\LL_\lambda$ the corresponding ladder,  and $r=r(n) \subset \LL_\lambda$ a qi  ray, then there exists
$z\in\partial \til{M^h}$ such that $r(n) \rightarrow z = \partial i(\lambda_{-\infty }) =
\partial i(\lambda_{\infty })$ as $n \rightarrow \infty$. 
\label{qgeodasymp-cusp}
\end{prop}

\begin{proof} We first observe that both end-points $\lambda_{-\infty }, \lambda_{\infty }$ of the CT leaf $\lambda$ cannot be parabolics.
For then they would have to be base points of {\it different} horoballs in $\til M$ as they correspond to different lifts of the cusp(s) of $M$. 

\smallskip

\noindent {\bf Case 1:} Both $\lambda_{-\infty }, \lambda_{\infty }$ are non-parabolic. \\
The proof  of Proposition \ref{qgeodasymp} goes through in this context mutatis mutandis.

\smallskip

\noindent {\bf Case 2:} Exactly one of $\lambda_{-\infty }, \lambda_{\infty }$  is a parabolic. \\
Without loss of generality assume that  $\lambda_{-\infty }$ is a parabolic. 
Let $B$ be the horoball in $\til{M^h}$ based at $w = \partial i (\lambda_{-\infty} )$ and let $H$ be the horosphere boundary of $B$. Let $o$ be
the point of intersection of $\lambda$ with $H$. For $p, q \in \til{ M^h}$, $(p,q)_h$ and $\overline{pq}$ will denote respectively
geodesics in $(\til{ M^h},d)$ and $(\til{ M^h},d_G)$.

Choose a sequence of points $a_n, b_n \in \lambda$
such that $a_n \rightarrow \lambda_{-\infty} $
and $b_n \rightarrow \lambda_\infty $.  Then by the existence
of Cannon-Thurston maps for $i: \til{S^h} \rightarrow \til{M^h}$ (Theorem \ref{split}) it follows  that there exists
a function $M(n) \rightarrow \infty$ as $n \rightarrow \infty$ such that $(a_n,b_n)_h$ lies outside $B_{M(n)} (o) \subset \til{M^h}$.
Hence, if $q_n = (a_n,b_n)_h \cap H$ then $d(q_n, o) \geq M(n)$ and the geodesic subsegment $(q_n,b_n)_h$ lies outside $B_{M(n)} (o) \subset \til{M^h}$.

Let $N = \til{M^h} \setminus \bigcup_\alpha B_\alpha$ be the complement of open horoballs and $d_G$ be the graph metric on $N$ obtained after first
partially electrocuting horospheres (cf. Section \ref{hypsec}). By Lemma \ref{ea-genl}
$(q_n,b_n)_h$ and $\overline{q_nb_n}$ lie in a bounded neighborhood of each other in $(N,d_G)$. 

The $d_G$-distance $d_G(o, q_n)$ is  equal to
the number of vertical
blocks between $o$ and $q_n$. But $a_n \rightarrow \lambda_{-\infty}$ implies
$q_n \rightarrow \infty$ in $\til{M^h}$. Hence $d_G(o, q_n) \rightarrow \infty$ as $a_n \rightarrow \lambda_{-\infty}$. 

By Corollary \ref{ctleafinf}, $d_G(o,b_n)\rightarrow \infty$ as $n\rightarrow \infty$. Hence by Proposition \ref{ctdg} there exists
a function $M_1(n) \rightarrow \infty$ as $n \rightarrow \infty$ such that $\overline{q_nb_n}$ lies outside $B_{M_1(n)} (o) \subset (N,d_G)$.

Now recall that $\Pi_\lambda : N \rightarrow \LL_\lambda$ is a coarse Lipschitz retract by Theorem \ref{retract}. Hence $\Pi_\lambda [\overline{q_nb_n}]
\subset \LL_\lambda$ is a uniform quasigeodesic in $(N,d_G)$.

Further, since $q_n$ belongs to $H$ and since $\Pi_\lambda$ essentially fixes the horosphere $H$, it follows that $d_G( \Pi_\lambda (q_n), q_n) \leq 1$.
Also $\Pi_\lambda (b_n) = b_n$. Therefore there exists
a function $M_2(n) \rightarrow \infty$ as $n \rightarrow \infty$ such that $\Pi_\lambda [\overline{q_nb_n}]$ lies outside $B_{M_2(n)} (o) \subset (N,d_G)$.

Next, since $H \cap \LL_\lambda$ and $b_n$ lie on different sides of the  qi  ray $r=r(n) \subset \LL_\lambda$ it follows that there exists
$z_n \in \overline{q_nb_n}$ such that $d_G(z_n,r)$ is uniformly bounded.

Also there exists $t_n \in (q_n,b_n)_h$ such that  $d_G(z_n,t_n)$ and hence $d_G(t_n,r)$ is uniformly bounded.

Since $t_n \in (q_n,b_n)_h$ it follows that $t_n \rightarrow \partial i (\lambda_{-\infty} ) = \partial i (\lambda_{\infty} )$.
Since $d_G(t_n,r)$ is uniformly bounded, there exists $s_n \in r$ such that
$d_G(t_n,s_n)$ is uniformly bounded and therefore $t_n, s_n$ are separated by a uniformly bounded number of split components. By uniform graph
quasiconvexity of split components (Theorem \ref{gqc}) it follows that $s_n \rightarrow \partial i (\lambda_{-\infty} ) = \partial i (\lambda_{\infty} )$. 

Finally if $r_{s_n}$ denotes the part of the ray $r$ `above' $s_n$, (i.e. $[s_n, \infty)$) then joining  points of $r_{s_n}$ in
successive blocks by hyperbolic geodesics we obtain an electro-ambient quasigeodesic $\sigma_n$. By Lemma \ref{ea-genl}
there exist hyperbolic geodesics $\tau_{m,n}$ joining $r(m), r(n)$ for $m > n$ and contained in
a bounded neighborhood of $\sigma_n \cup B$ in $\til{M^h}$. Hence  $r(n) \rightarrow  \partial i(\lambda_{-\infty }) =
\partial i(\lambda_{\infty })$ as $n \rightarrow \infty$.
\end{proof}

We are now in a position to prove the analogue of Theorem \ref{main} for surfaces with cusps.

\begin{theorem} \label{main-cusp} 
 Let $\partial i(a) = \partial i(b)$ for $a, b \in
S^1_\infty$ be two distinct points that are identified by the
Cannon-Thurston map corresponding to a simply degenerate  surface
group (without accidental parabolics). 
Then $a, b$ are either ideal 
end-points of a leaf of the ending lamination, 
or ideal boundary points of a
complementary ideal polygon. Further, if $a, b$ are either ideal 
end-points of a leaf of a lamination, or ideal boundary points of a
complementary ideal polygon, then $\partial i(a) = \partial i(b)$. \end{theorem}

\begin{proof}  The second statement has been shown in Proposition \ref{ptpreimage}. \\ To prove the first statement,
it suffices to show that $\RCT$ is unlinked.
Suppose now that $\lambda$ and $\mu$ are intersecting CT leaves, i.e.
$\partial i(\lambda_{-\infty }) =
\partial i(\lambda_{\infty })$  and $\partial i(\mu_{-\infty }) =
\partial i(\mu_{\infty })$. 

Consider the ladders ${\LL}_\lambda$ and ${\LL}_\mu$. Let $r(i) \subset \lambda_i \cap \mu_i$ be a 
quasigeodesic ray as per Remark \ref{twin-cusp}. 
By Proposition \ref{qgeodasymp-cusp}, $r$ converges
to a point $z$ on $\partial \til{M}$ 
such that $z =  \partial i(\lambda_{-\infty }) =
\partial i(\lambda_{\infty })  = \partial i(\mu_{-\infty }) =      
\partial i(\mu_{\infty })$. Hence if $ \{ a, b \},  \{ c, d \} \in \RCT$, then either $ \{ a, b, c, d \} $ are all mutually related in $\RCT$, or 
$ \{ a, b \},  \{ c, d \}$ are unlinked. By Lemma \ref{unlink}, $\RCT$ is induced by a lamination $\Lambda_{CT}$. By Proposition \ref{ptpreimage}, 
the ending lamination $\Lambda_{EL}$ is contained in $\Lambda_{CT}$. Since $\Lambda_{EL}$ is filling and arational, it follows that
$\Lambda_{EL} = \Lambda_{CT}$.
\end{proof}

The 
modifications necessary to pass from the simply degenerate case to the totally degenerate case are exactly as in the case of surfaces without cusps.

\bibliography{ptpre}
\bibliographystyle{alpha}

\end{document}